\documentclass{proc-l}
\usepackage{graphicx}
\newtheorem{theorem}{Theorem}[section]
\newtheorem{lemma}[theorem]{Lemma}
\newtheorem{corollary}[theorem]{Corollary}

\theoremstyle{definition}
\newtheorem{definition}[theorem]{Definition}
\newtheorem{example}[theorem]{Example}

\theoremstyle{remark}

\numberwithin{equation}{section}

\newcommand{\floor}[1]{\mbox{$\lfloor{#1}\rfloor$}}

\begin{document}

\title{The inverse rook problem on Ferrers boards}

\author{Abigail G. Mitchell}
\address{Department of Mathematics\\
University of Notre Dame\\
Notre Dame, Indiana  46556-5683}
\email{mitchell.85@nd.edu}

\subjclass[2000]{<subjectclass>}
\date{\today}



\begin{abstract}
Rook polynomials have been studied extensively since 1946, principally as a method for enumerating restricted permutations.  However, they have also been shown to have many fruitful connections with other areas of mathematics, including graph theory, hypergeometric series, and algebraic geometry.  It is known that the rook polynomial of any board can be computed recursively. \cite{ri58,mi04b}

The naturally arising inverse question --- given a polynomial, what board (if any) is associated with it? --- remains open.  In this paper, we solve the inverse problem completely for the class of Ferrers boards, and show that the increasing Ferrers board constructed from a polynomial is unique.

\end{abstract}

\maketitle

\section{Introduction}

Rook polynomials provide a method of enumerating permutations with
restricted position.  Their study was begun in 1946 by Kaplansky
and Riordan~\cite{ka46r}, with applications to card-matching problems.
Riordan's 1958 book~\cite{ri58} is considered the first systematic
analysis, and remains a classic treatment of the subject.  A series
of papers by Goldman et al. (\cite{go75jw,go76jrw,go78jw,go77jw,go76jw}) in the 1970s expanded the field by
applying more advanced combinatorial methods.  

In this paper, we shall wish to restrict our attention to the class of Ferrers boards.  Since 1975 (\cite{go75jw}) these have occupied a prominent place in the literature.  This is principally due to the depth and variety of the connections Ferrers boards exhibit to other parts of mathematics.  Ferrers boards are related to chromatic polynomials \cite{br94rw}, algebraic geometry \cite{di97}, hypergeometric series \cite{ha96a}, permutation statistics \cite{bu04}, quantum mechanical operators \cite{va04}, and several types of digraph polynomials \cite{ch95g,da00m,fa91w}.  Such diverse applications make the inverse problem an especially pressing question for Ferrers boards.

In this paper, we solve the inverse problem completely for Ferrers boards.  Our proof is constructive, and using the Foata-Sch\"utzenberger Theorem~\cite{fo70s}, we show that the board thus constructed is unique within the subclass of increasing Ferrers boards.

\section{Preliminaries}\label{sec:prelim}

Let $B$ be a generalized chessboard -- a set of cells arranged in rows
and columns.  The rook is a chess piece which attacks on rows and
columns; by a rook placement we shall mean a non-attacking placement
of $k$ indistinguishable rooks on the board $B$.

This intuitive definition can be formalized in several ways.
A board $B$ may be regarded as a subset of $[1,2,\ldots ,m]\times [1,2,\ldots ,n]$
for some $m,n \in \mathbb{N} $.
  A rook placement on $B$ then corresponds to a choice of $k$ elements
of $B$, ${(x_1,y_1),(x_2,y_2),\ldots ,(x_k,y_k)}$ such that ${x_i}={x_j}$ or
${y_i}={y_j}$ implies $i=j$.  

Alternatively, let $D\subseteq [1,2,\ldots ,m]$ be an arbitrary subset.  Then a rook placement is an injective function
 $f:D \rightarrow [1,2,\ldots ,n]$
 such that, for any $i\in D$, $(i,f(i))\in B$.

$B$ may also be seen as an element of ${M_{m,n}}[\mathbb{F}_2]$, an
$m \times n$ matrix with binary entries.  We write in this case $B=(b_{i,j})$.  In this interpretation, a
placement of $k$ rooks on $B$ corresponds to a choice of $k$
independent 1's in $B$.

A third formalization is to regard $B$ as a bipartite graph on vertex
sets $[1,2,...,m]$ and $[1,2,...,n]$, where the graph contains edge
$(i,j)$ iff $(i,j)$ is a cell in the corresponding chessboard.  A
rook placement then corresponds to a partial matching on the graph.

In this paper, we will not often need to invoke any of these formalizations.  Insofar as we use any of them, we will think of a board as an $m\times n$ binary matrix; this serves to ground the intuitive concepts of `row' and `column', as well as providing a convenient notation for identifying individual cells.

\begin{definition}
The (classical) rook polynomial of a board $B$ is the ordinary generating function
\[
R(B;x)=\sum_{k=0}^{\infty}r_k x^k,
\]
where the coefficient $r_k$ is the number of placements of $k$ rooks
on $B$.  When its omission will cause no confusion, we will assume
the variable to be $x$, and will simply write $R(B)$.
\end{definition}

\begin{definition}(after Goldman et al. in \cite{go76jw})
For $n\in \mathbb{N}$, the $n$-factorial rook polynomial of a board $B$ is defined as
\[
p_n(B;x)=\sum_{k=0}^n r_k (x)_{n-k},
\]
where $(x)_{i}$ is the falling factorial $x(x-1)(x-2)\ldots (x-i+1)$.  As with the classical rook polynomial, we shall usually assume the variable $x$ and write $p_n(B)$.
\end{definition}

A Ferrers board is a board made up of a sequence of adjacent solid columns of nondecreasing height, with a common lower edge.  (We shall find it convenient to allow columns of zero height.)  If the sequence of column heights is strictly increasing, the board is called an increasing Ferrers board.  

As a binary matrix, $B=(b_{i,j})$ is a Ferrers board if there exists a nondecreasing sequence of nonnegative integers $h_1,h_2,\ldots ,h_n$ such that $b_{i,j} = 1$ if and only if $j \leq h_i$.

\begin{figure}
\begin{center}
  \includegraphics[totalheight=0.5in]{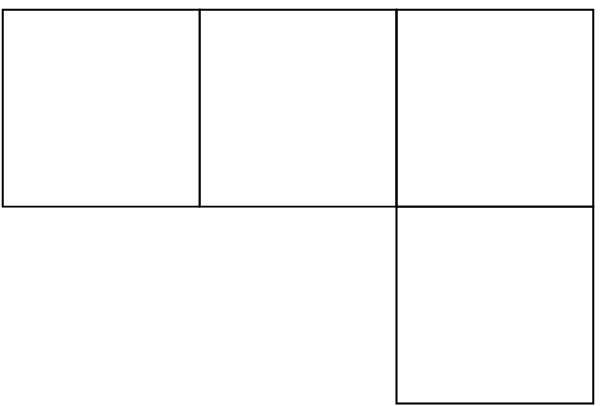}
  \hspace{1.0in}
  \includegraphics[totalheight=0.5in]{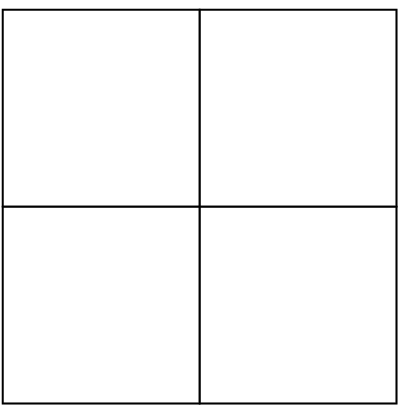}
\caption{Two rook equivalent boards, with rook polynomial $1+4x+2x^2$.} \label{equivalent-boards}
\end{center}
\end{figure}

Two boards $A$,$B$ are said to be rook equivalent if $R(A)=R(B)$.  A
sufficient condition for rook equivalence is that $B$ can be obtained
from $A$ by permutation of rows and columns.  (This condition is not
necessary, as demonstrated by the two boards in Figure~\ref{equivalent-boards}.)

Several immediate observations can be made regarding the coefficients $r_k$ of a rook polynomial.
\begin{theorem}
Let $B$ be a board, and let $R(B;x)=\sum _{k=0} ^\infty r_k x^k$ be its rook polynomial.
\begin{enumerate}
\item $r_k \geq 0$ for all $k$.
\item $r_k=0$ implies $r_j=0$ for all $j \geq k$.
\item $r_0=1$.
\item $r_1=|B|$, the number of cells in $B$.
\item $r_k \leq {r_1\choose k}$ for all $k$.
\item $r_k \leq {r_{k-1}\choose 2}$ for all $k \geq 2$.
\end{enumerate}
\end{theorem}
\begin{proof}
Most of these are immediate.  For the last, note that any placement of $k$ rooks gives rise to $k$ placements of $k-1$ rooks, by removing each rook in turn.  For $k\geq 2$, any pair of such ($k-1$)-placements is sufficient to uniquely determine a $k$-placement.  An upper bound for the number of placements of $k$ rooks is therefore ${r_{k-1}\choose 2}$.
\end{proof}

These conditions are necessary, but not sufficient, for a given polynomial to be the rook polynomial of some board.  For instance, $1+4x+x^2$ and $1+4x+5x^2+x^3$ satisfy all the above conditions, but are not rook polynomials.  In the next section, we establish conditions which are both necessary and sufficient for a polynomial to be the rook polynomial of a Ferrers board.

\section{Results}\label{sec:main}

We shall need the following definitions and results:

\begin{definition} (after \cite{go75jw})
Let $B$ be a Ferrers board with $c$ columns, of heights $h_1,h_2,\ldots ,h_c$.
We define the height vector $h(B)$ to be $h(B)=(h_1,h_2,\ldots ,h_c)$.
For $n\geq c$ we define the $n$-height vector $h_n(B)$ to be
$h_n(B)=(h_1^{(n)},h_2^{(n)},\ldots ,h_c^{(n)})$, where $h_i^{(n)}=0$ for
$i=1,2,\ldots ,n-c$ and $h_i^{(n)}=h_{i-(n-c)}$ for $i=n-c+1,\ldots ,n$.
We also define the $n$-structure vector $s_n(B)=(s_1^{(n)},s_2^{(n)},\ldots ,s_n^{(n)})$ where $s_i^{(n)}=h_i^{(n)}-(i-1)$.
\end{definition}

\begin{example}
Consider the following Ferrers board.
\begin{center}
\raisebox{0.45in} {$B\; :=\;$}
\includegraphics[totalheight=1.0in]{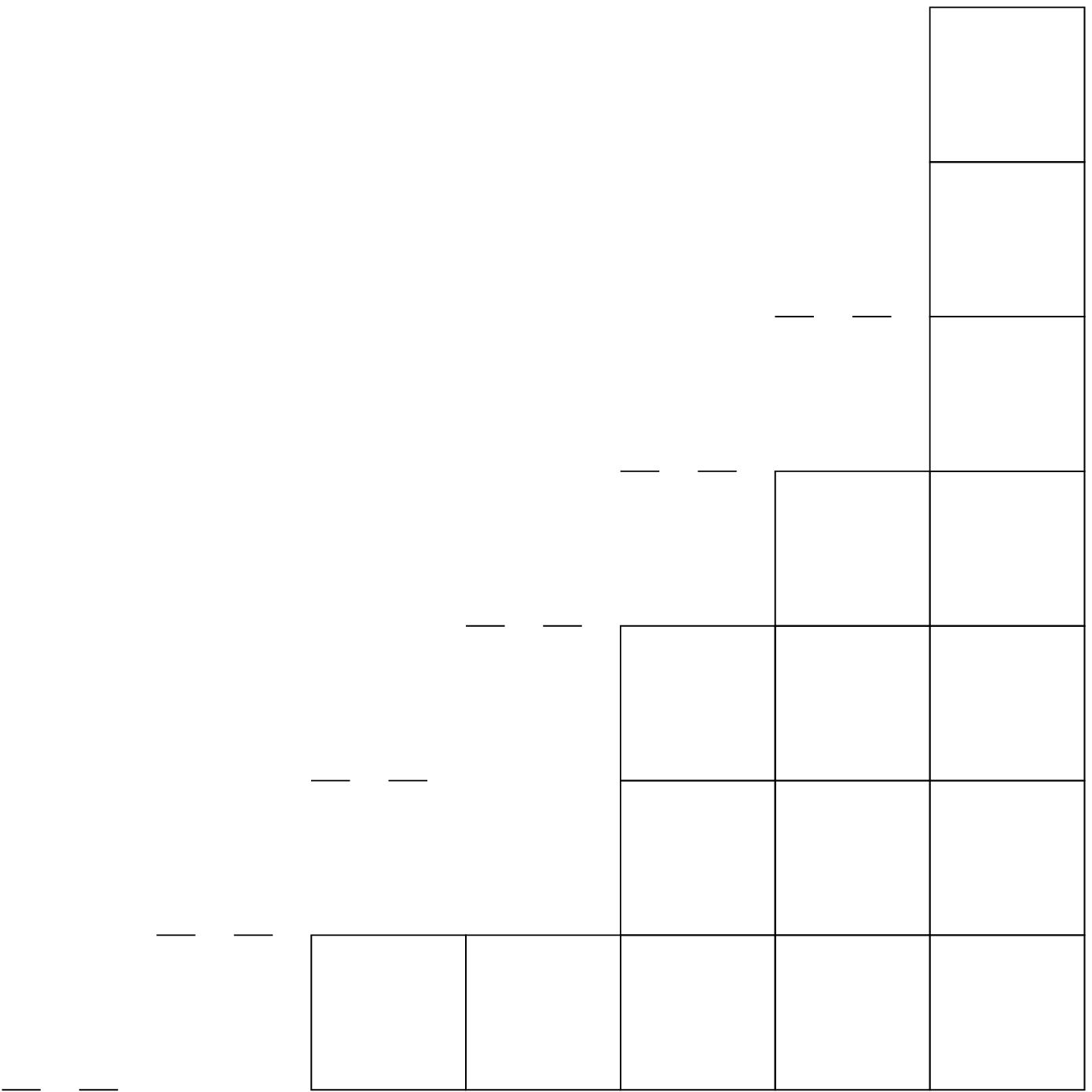}
\end{center}
The height vector of $B$ is $h(B)=(1,1,3,4,7)$.  Its $7$-height vector is $h_7(B)=(0,0,1,1,3,4,7)$, and its $7$-structure vector is $s_7(B)=(0,-1,-1,-2,-1,-1,1)$.

\end{example}

Note that, as shown in this example, $h_n(B)$ corresponds to placing a set of $n-c$ columns of zero height to the left of $B$, and that the components of $s_n(B)$ correspond to heights above or below the diagonal of this augmented board (the dotted line shown).

\begin{theorem}[Factorization Theorem. Goldman et al., Theorem 2 in \cite{go75jw}]
\label{factorization}
Let $B$ be a Ferrers board with $c$ columns.  For $n\geq c$, let
$p_n(B;x)$ be the $n$-factorial rook polynomial of $B$, and
 $(s_1^{(n)},s_2^{(n)},\ldots ,s_n^{(n)})$
 be the $n$-structure vector of $B$.  Then
\[
p_n(B;x)=\prod_{i=1}^n (x+s_i^{(n)}).
\]
\end{theorem}

\begin{corollary} [Goldman et al. \cite{go75jw}]
Two Ferrers boards are rook equivalent if and only if, for some $n$, their $n$-factorial rook polynomials are equal.
\end{corollary}

\begin{theorem}[Foata-Sch\"utzenberger \cite{fo70s}]
\label{increasing}
Every Ferrers board is rook equivalent to a unique increasing Ferrers board.
\end{theorem}

\begin{definition}
For any $q(x) = \sum_{k=0}^{m} q_k x^k \: \in \mathbb{Z} [x]$, and for $n\geq m$, we associate to $q$ the $n$-factorial polynomial $q_n$ defined as
\[
q_n(x)=\sum_{k=0}^{m} q_k (x)_{n-k},
\]
where $(x)_{i}$ denotes the falling factorial.
\end{definition}

\begin{figure}
\begin{center}
  \includegraphics[totalheight=1.0in]{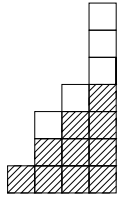}
  \caption{An increasing Ferrers board with 4 columns.} \label{triangle}
\end{center}
\end{figure}

\begin{lemma}\label{columns}
Let $B$ be an increasing Ferrers board having $n$ cells.  Then the number of nonzero columns of $B$ is strictly less than $\sqrt{2n}$.
\end{lemma}
\begin{proof}
Let $c$ be the number of columns of $B$.  Since $B$ is strictly increasing, $B$ contains cells forming a $c \times c$ right triangle, as shown in Figure~\ref{triangle}.  The number of cells in this triangle is $frac{c(c+1)}{2}$, so we have
\begin{align*}
\frac{c(c+1)}{2} &\leq n \\
c^2+c-2n &\leq 0 \\
\frac{-1-\sqrt{1+8n}}{2} &\leq c \leq \frac{-1+\sqrt{1+8n}}{2} \\
\end{align*}
Since $c\geq 0$, this gives
\[
c\leq \frac{-1+\sqrt{1+8n}}{2} < \sqrt{2n}.
\]
\end{proof}

Now we are ready to state the main theorem of this section.

\begin{theorem}
Let $q(x)=\sum_{k=0}^m q_k x^k$ be a polynomial over $\mathbb{Z}$.  
Let $n=\floor{\sqrt{2q_1}}$, and let $p(x)=\sum_{k=0}^m q_k (x)_{n-k}$ be the n-factorial polynomial associated with $q(x)$.  $q(x)$ is the rook polynomial of an increasing Ferrers board $B$ if and only if the following conditions are satisfied:
\begin{enumerate}
\item $p(x)$ has only integer roots.
\item $(x)_{t+1}$ divides $p(x)$, where $t$ is the largest root of $p(x)$.
\end{enumerate}
Furthermore, if it exists $B$ is unique, and has column height vector $h_n = (-x_1+t+1,-x_2+t+2,\ldots ,-x_c+t+c)$, where $x_i$ are the (possibly repeated) roots of $\frac{p(x)}{(x)_{t+1}}$, arranged in nonincreasing order.
\end{theorem}

\begin{proof}
Assume such a board $B$ exists.  Then it must have $q_1$ cells, so by Lemma~\ref{columns}
it has no more than $n$ columns.
The $n$-structure vector of any increasing Ferrers board is of the form $(0,-1,-2,\ldots,-t,u_1,u_2,\ldots,u_c)$ with $-t \leq u_1 \leq u_2 \leq \cdots \leq u_c$. 
By the hypothesis that $q(x)=R(B;x)$, $p(x)$ is the $n$-factorial rook polynomial of $B$.  Therefore, by Theorem~\ref{factorization},
\begin{align*}
p(x) &= \prod_{i=1}^{n} (x+s_i^{(n)})\\
     &= x(x-1)(x-2)\cdots (x-t)(x-u_1)(x-u_2)\cdots (x-u_c)\\
     &= (x)_{t+1} (x-u_1)(x-u_2)\cdots (x-u_c).\\
\end{align*}
Thus, $(x)_{t+1}$ divides $p(x)$, and all the roots of $p(x)$ are integers.  By Theorem~\ref{increasing}, $B$ is unique, and it follows from the definitions that $B$ has the column heights stated.

For the converse, assume $p(x)$ satisfies the given conditions.  Then
\[
p(x) = (x)_{t+1} \hat{p}(x),
\]
where $\hat{p}(x)$ is a polynomial of degree $n-(t+1)$ with integer roots.  Let $u_1,u_2,\ldots u_{n-(t+1)}$ be its roots, arranged in nonincreasing order.  Then
\begin{align*}
p(x) &= (x)_{t+1} \prod_{i=1}^{n-(t+1)} (x-u_i)\\
     &= \prod_{i=1}^{t+1} (x-(i-1)) \prod_{i=1}^{n-(t+1)} (x-u_i)\\
     &= \prod_{i=1}^n (x+s_i),\\
\end{align*}
where $s_i=-(i-1)$ for $i = 1,2,\ldots ,t+1$ and $s_i=-u_i$ for $i = t+2,\ldots ,n$.
Let $s_n=(s_1,s_2,\ldots ,s_n)$.  This is exactly the $n$-structure vector of an increasing Ferrers board $B$ with column height vector $h_n = (-u_1+t+1,-u_2+t+2,\ldots ,-u_c+t+c)$.  Therefore, by Theorem~\ref{factorization}, $p(x)$ is the $n$-factorial rook polynomial of $B$, and so $q(x)$ is the classical rook polynomial of $B$.  By Theorem~\ref{increasing}, $B$ is uniquely determined.  
\end{proof}

The following corollary is immediate, by Theorem~\ref{increasing}.
\begin{corollary}
The same conditions are both necessary and sufficient for $q(x)$ to be the rook polynomial of {\em any} Ferrers board.
\end{corollary}

\section{Conclusions} \label{conclusions}
The ascendancy of Ferrers boards in the literature can be attributed to various factors.  One is the real frequency with which Ferrers board problems crop up in other fields; another is that more direct combinatorial interpretations (for example, as partitions) make their study more tractable.  Still another factor is the apparent tendency of new mathematics to accumulate around existing results.

This paper falls primarily into the latter two categories.  The enduring openness of the inverse problem motivated this research, and the combined results of Foata-Sch\"utzenberger and Goldman-Joichi-White presented what seems a natural approach to solving it, at least for Ferrers boards.  Unfortunately, it does not seem that the results obtained in this paper will be readily extensible to more general boards.

\bibliographystyle{amsplain}

\end{document}